\documentclass[12pt,a4paper]{article}

\usepackage{amsmath}
\usepackage{amssymb}
\usepackage{mathrsfs}
\usepackage{bbold}
\usepackage{graphicx}

\renewcommand{\leq}{\leqslant}

\renewcommand{\geq}{\geqslant}

\makeatletter
\newcommand*{\bigcorr@macro}[2]{\sbox{0}{\mbox{$#1($}}\dimen0=\ht0
                \advance\dimen0 by \dp0
                \multiply\dimen0 by #2 \divide\dimen0 by 100}
\newcommand*{\bigcorr@big}[2]{\mbox{$#1\left#2\bigcorr@macro{#1}{85}\vrule
                   height \dimen0 depth 0pt width 0pt\right.\n@space$}}
\newcommand*{\bigcorr@Big}[2]{\mbox{$#1\left#2\bigcorr@macro{#1}{115}\vrule
                   height \dimen0 depth 0pt width 0pt\right.\n@space$}}
\newcommand*{\bigcorr@bigg}[2]{\mbox{$#1\left#2\bigcorr@macro{#1}{145}\vrule
                   height \dimen0 depth 0pt width 0pt\right.\n@space$}}
\newcommand*{\bigcorr@Bigg}[2]{\mbox{$#1\left#2\bigcorr@macro{#1}{175}\vrule
                   height \dimen0 depth 0pt width 0pt\right.\n@space$}}
\DeclareRobustCommand*{\big}[1]{{\mathpalette\bigcorr@big{#1}}}
\DeclareRobustCommand*{\Big}[1]{{\mathpalette\bigcorr@Big{#1}}}
\DeclareRobustCommand*{\bigg}[1]{{\mathpalette\bigcorr@bigg{#1}}}
\DeclareRobustCommand*{\Bigg}[1]{{\mathpalette\bigcorr@Bigg{#1}}}
\DeclareRobustCommand*{\bi}[1]{{#1}}
\DeclareRobustCommand*{\bil}[1]{\mathopen{\bi{#1}}}

\DeclareRobustCommand*{\bir}[1]{\mathclose{\bi{#1}}}
\makeatother

\newtheorem{theorem}{Theorem}
\newtheorem{lemma}[theorem]{Lemma}

\newcommand{\cN}{\mathcal{N}}
\newcommand{\sF}{\mathscr{F}}
\newcommand{\RR}{\mathbb{R}}
\newcommand{\NN}{\mathbb{N}}
\newcommand{\eb}{\mathbf{E}}
\newcommand{\pb}{\mathbf{P}}
\newcommand{\1}{\mathbb{1}}
\newcommand{\argsup}{\mathop{\rm argsup}\limits}
\newcommand{\argmax}{\mathop{\rm argmax}\limits}
\newcommand*{\abs}[1]{\left|#1\right|}
\newcommand*{\mybiblink}[1]{{\par\vglue -3pt\scriptsize\texttt{#1}\par}}

\begin{document}

\title{On Compound Poisson Processes Arising in\\
       Change-Point Type Statistical Models\\
       as Limiting Likelihood Ratios}
\author{Sergue\"{\i} \textsc{Dachian}\thanks{Laboratoire de Math\'ematiques UMR6620,
Universit\'e Blaise Pascal, Clermont Universit\'e, F--63177 Aubi\`ere CEDEX. E-mail:
\texttt{Serguei.Dachian@math.univ-bpclermont.fr}}
\and
Ilia \textsc{Negri}\thanks{Department of Information Technology and
Mathematical Methods, University of Bergamo, I--24044 Dalmine (BG). E-mail:
\texttt{ilia.negri@unibg.it}}}

\date{}
\maketitle

\begin{abstract}
Different change-point type models encountered in statistical inference for
stochastic processes give rise to different limiting likelihood ratio
processes. In a previous paper of one of the authors it was established that
one of these likelihood ratios, which is an exponential functional of a
two-sided Poisson process driven by some parameter, can be approximated (for
sufficiently small values of the parameter) by another one, which is an
exponential functional of a two-sided Brownian motion. In this paper we
consider yet another likelihood ratio, which is the exponent of a two-sided
compound Poisson process driven by some parameter. We establish, that
similarly to the Poisson type one, the compound Poisson type likelihood ratio
can be approximated by the Brownian type one for sufficiently small values of
the parameter. We equally discuss the asymptotics for large values of the
parameter and illustrate the results by numerical simulations.

\end{abstract}



\section{Introduction}

In this work we are interested by the asymptotic study of non-regular
parametric statistical models encountered in statistical inference for
stochastic processes. An exhaustive exposition of the parameter estimation
theory in both regular and non-regular cases is given in the classical
book~\cite{IKh} by Ibragimov and Khasminskii. They have developed a general
theory of estimation based on the analysis of renormalized likelihood
ratio. Their approach consists in proving first that the renormalized
likelihood ratio (with a properly chosen renormalization rate) weekly
converges to some non-degenerate limit: the limiting likelihood ratio
process. Thereafter, the properties of the estimators (namely their rate of
convergence and limiting distributions) are deduced. Finally, based on the
estimators, one can also construct confidence intervals, tests, and so on.
Note that this approach also provides the convergence of moments, allowing one
to deduce equally the asymptotics of some statistically important quantities,
such as the mean squared errors of the estimators.

It is well known that in the regular case the limiting likelihood ratio is
given by the LAN property and is the same for different models (the
renormalization rate being usually $1/\sqrt{n}\,$). So, the classical
estimators~--- the maximum likelihood estimator and the Bayesian
estimators~--- are consistent, asymptotically normal (usually with rate
$1/\sqrt{n}\,$) and asymptotically efficient.

In non-regular cases the situation essentially changes: the renormalization
rate is usually better (for example, $1/n$ in change-point type models), but
the limiting likelihood ratio can be different in different models. So, the
classical estimators are still consistent, but may have different limiting
distributions (though with a better rate) and, in general, only the Bayesian
estimators are asymptotically efficient.

In~\cite{D} a relation between two different limiting likelihood ratios
arising in change-point type models was established by one of the
authors. More precisely, it was shown that the first one, which is an
exponential functional of a two-sided Poisson process driven by some
parameter, can be approximated (for sufficiently small values of the
parameter) by the second one, defined by
\begin{equation}
\label{proc2}
Z_0(x)=\exp\left\{W(x)-\frac12\abs x\right\},\quad x\in\RR,
\end{equation}
where $W$ is a standard two-sided Brownian motion. In this paper we consider
yet another limiting likelihood ratio process arising in change-point type
models and show that it is related to $Z_0$ in a similar way.

\subsection*{The process \boldmath$Z_{\gamma,f}$}

We introduce the random process $Z_{\gamma,f}$ on $\RR$ as the exponent of a
two-sided compound Poisson process given by
\begin{equation}
\label{proc1}
\ln Z_{\gamma,f}(x)=\begin{cases}
\vphantom{\bigg)}\sum_{k=1}^{\Pi_+(x)}\ln\frac{f(\varepsilon_k^++\gamma)}
{f(\varepsilon_k^+)}\,,&\text{if } x\geq 0,\\
\vphantom{\bigg)}\sum_{k=1}^{\Pi_-(-x)}\ln\frac{f(\varepsilon_k^--\gamma)}
{f(\varepsilon_k^-)}\,,&\text{if } x\leq 0,\\
\end{cases}
\end{equation}
where $\gamma>0$, $f$ is a strictly positive density of some random variable
$\varepsilon$ with mean $0$ and variance $1$, $\Pi_+$ and $\Pi_-$ are two
independent Poisson processes of intensity $1$ on $\RR_+$,
$\varepsilon_k^{\pm}$ are independent random variables with density $f$ which
are also independent of $\Pi_\pm$, and we use the convention
$\sum_{k=1}^{0}a_k=0$. We equally introduce the random variables
\begin{equation}
\label{vars1}
\begin{aligned}
\zeta_{\gamma,f}&=\frac{\int_{\RR}x\,Z_{\gamma,f}(x)\;dx}
{\int_{\RR}\,Z_{\gamma,f}(x)\;dx}\,,\\
\xi_{\gamma,f}^-&=\inf\Bigl\{z :
Z_{\gamma,f}(z)=\sup_{x\in\RR}Z_{\gamma,f}(x)\Bigr\},\\
\xi_{\gamma,f}^+&=\sup\Bigl\{z :
Z_{\gamma,f}(z)=\sup_{x\in\RR}Z_{\gamma,f}(x)\Bigr\},\\
\xi_{\gamma,f}^\alpha&=\alpha\,\xi_{\gamma,f}^-+
(1-\alpha)\,\xi_{\gamma,f}^+,\quad\alpha\in[0,1],
\end{aligned}
\end{equation}
related to this process, as well as their second moments
$B_{\gamma,f}=\eb\zeta_{\gamma,f}^2$ and
$M_{\gamma,f}^\alpha=\eb(\xi_{\gamma,f}^\alpha)^2$.

An important particular case of this process is the one where the density $f$
is Gaussian, that is, $\varepsilon\sim\cN(0,1)$. In this case we will omit the
index $f$ and write $Z_\gamma$ instead of $Z_{\gamma,f}$, $\xi_\gamma^\alpha$
instead of $\xi_{\gamma,f}^\alpha$, and so on. Note that since
$$
\ln\frac{f(\varepsilon\pm\gamma)}{f(\varepsilon)}=
\mp\gamma\varepsilon-\frac{\gamma^2}{2}\sim\cN(-\gamma^2/2,\gamma^2),
$$
the process $Z_\gamma$ is symmetric and has Gaussian jumps.

The process $Z_{\gamma,f}$, up to a linear time change, arises in some
non-regular, namely change-point type, statistical models as the limiting
likelihood ratio process, and the variables $\zeta_{\gamma,f}$ and
$\xi_{\gamma,f}^\alpha$ as the limiting distributions of the Bayesian
estimators and of the appropriately chosen maximum likelihood estimator,
respectively. The maximum likelihood estimator being not unique in the
underlying models, the appropriate choice here is a linear combination with
weights $\alpha$ and $1-\alpha$ of its minimal and maximal values. Moreover,
the quantities $B_{\gamma,f}$ and $M_{\gamma,f}^\alpha$ are the limiting mean
squared errors (sometimes also called limiting variances) of these estimators
and, the Bayesian estimators being asymptotically efficient, the ratio
$E_{\gamma,f}^\alpha=B_{\gamma,f}/M_{\gamma,f}^\alpha$ is the asymptotic
relative efficiency of this maximum likelihood estimator.

The examples include the two-phase regression model and the threshold
autoregressive (TAR) model. The linear case of the former was studied by Koul
and Qian in~\cite{KQ}, while the non-linear one was investigated by Ciuperca
in~\cite{Ciu}. Concerning the TAR model, the first results were obtained by
K.S.~Chan in~\cite{Cha}, while a more recent study was performed by N.H.~Chan
and Kutoyants in~\cite{CK}. Note however, that the estimator studied
in~\cite{Cha} is the least squares estimator (which is, in the Gaussian case,
equivalent to the maximum likelihood estimator), while the model considered
in~\cite{CK} is the Gaussian TAR model. So, only the processes~$Z_\gamma$ are
known to arise as limiting likelihood ratios in the TAR model.  Note also that
in both models, the parameter~$\gamma$ of the limiting likelihood ratio is
related to the jump size of the model.

\subsection*{The process \boldmath$Z_0$}

On the other hand, many change-point type statistical models encountered in
various fields of statistical inference for stochastic processes rather have
as limiting likelihood ratio process, up to a linear time change, the process
$Z_0$ defined by~\eqref{proc2}. In this case, the limiting distributions of
the Bayesian estimators and of the maximum likelihood estimator are given by
\begin{equation}
\label{vars2}
\zeta_0=\frac{\int_{\RR}x\,Z_0(x)\;dx}{\int_{\RR}\,Z_0(x)\;dx}
\quad\text{and}\quad\xi_0=\argsup_{x\in\RR}Z_0(x),
\end{equation}
respectively, while the limiting mean squared errors of these estimators are
$B_0=\eb\zeta_0^2$ and $M_0=\eb\xi_0^2$. The Bayesian estimators are still
asymptotically efficient, and the asymptotic relative efficiency of the
maximum likelihood estimator is $E_0=B_0/M_0$.

A well-known example is the model of a discontinuous signal in a white
Gaussian noise exhaustively studied by Ibragimov and Khasminskii
in~\cite{IKh5} and~\cite[Chapter~7.2]{IKh}, but one can also cite change-point
type models of dynamical systems with small noise considered by Kutoyants
in~\cite{Kut2} and~\cite[Chapter~5]{Kut3}, those of ergodic diffusion
processes examined by Kutoyants in \cite[Chapter~3]{Kut4}, a change-point type
model of delay equations analyzed by K\"uchler and Kutoyants in~\cite{KK}, a
model of a discontinuous periodic signal in a time inhomogeneous diffusion
investigated by H\"opfner and Kutoyants in~\cite{HK}, and so on.

Let us also note that Terent'yev in~\cite{Ter} determined the Laplace
transform of $\pb\bigl(\abs{\xi_0}>t\bigr)$ and calculated the constant
$M_0=26$. Moreover, the explicit expression of the density of $\xi_0$ was
later successively provided by Bhattacharya and Brockwell in~\cite{BB}, by Yao
in~\cite{Yao} and by Fujii in~\cite{Fuj}. Regarding the constant $B_0$,
Ibragimov and Khasminskii in~\cite[Chapter~7.3]{IKh} showed by means of
numerical simulation that $B_0=19.5\pm0.5$, and so $E_0=0.73\pm0.03$. Later
in~\cite{Gol}, Golubev expressed $B_0$ in terms of the second derivative (with
respect to a parameter) of an improper integral of a composite function of
modified Hankel and Bessel functions. Finally in~\cite{RS}, Rubin and Song
obtained the exact values $B_0=16\,\zeta(3)$ and $E_0=8\,\zeta(3)/13$,
where~$\zeta$ is Riemann's zeta function defined by
$\zeta(s)=\sum_{n=1}^{\infty}1/n^s$.

\subsection*{The results of the present paper}

In this paper we establish that the limiting likelihood ratio processes
$Z_{\gamma,f}$ and $Z_0$ are related. More precisely, under some regularity
assumptions on $f$, we show that as $\gamma\to 0$, the process
$Z_{\gamma,f}(y/I\gamma^2)$, $y\in\RR$, (where $I$ is the Fisher information
related to $f$) converges weakly in the space $\mathcal{D}_0(-\infty
,+\infty)$ (the Skorohod space of functions on $\RR$ without discontinuities
of the second kind and vanishing at infinity) to the process~$Z_0$. Hence, the
random variables $I\gamma^2\zeta_{\gamma,f}$ and
$I\gamma^2\xi_{\gamma,f}^\alpha$ converge weakly to the random variables
$\zeta_0$ and $\xi_0$, respectively. We show equally that the convergence of
moments of these random variables holds and so, in particular, $I^2\gamma^4
B_{\gamma,f} \to 16\,\zeta(3)$, $I^2\gamma^4 M_{\gamma,f}^\alpha \to 26$ and
$E_{\gamma,f}^\alpha \to 8\,\zeta(3)/13$. Besides their theoretical interest,
these results have also some practical implications. For example, they allow
to construct tests and confidence intervals on the base of the distributions
of $\zeta_0$ and $\xi_0$ (rather than on the base of those of
$\zeta_{\gamma,f}$ and $\xi_{\gamma,f}^\alpha$, which depend on the density
$f$ and are not known explicitly) in models having the process $Z_{\gamma,f}$
with a small $\gamma$ as a limiting likelihood ratio. Also, the limiting mean
squared errors of the estimators and the asymptotic relative efficiency of the
maximum likelihood estimator can be approximated as
$$
B_{\gamma,f}\approx\frac{16\,\zeta(3)}{I^2\gamma^4}\,,\quad
M_{\gamma,f}^\alpha\approx\frac{26}{I^2\gamma^4} \quad\text{and}\quad
E_{\gamma,f}^\alpha\approx\frac{8\,\zeta(3)}{13}
$$
in such models.

These are the main results of the present paper, and they are presented in
Section~\ref{MR}, where we also briefly discuss the second possible
asymptotics $\gamma\to+\infty$ and present some numerical simulations of the
quantities $B_\gamma$, $M_\gamma^\alpha$ and $E_\gamma^\alpha$ for
$\gamma\in\left]0,\infty\right[$. Finally, the proofs of the necessary lemmas
are carried out in Section~\ref{PL}.

Concluding the introduction let us note that a preliminary exposition (in the
particular Gaussian case) of the results of the present paper can be found
in~\cite{DN} and~\cite{DN2}.

\section{Asymptotics of \boldmath$Z_{\gamma,f}$}
\label{MR}

Let $\gamma>0$, and let $f$ be a strictly positive density of some random
variable $\varepsilon$ with mean $0$ and variance $1$.

\subsection*{Regularity assumptions}

We will always suppose that $\sqrt{f}$ is continuously differentiable in
$L^2$, that is, there exists $\psi\in L^2$ satisfying
$\int_{\RR}\bigl(\sqrt{f(x+h)}-\sqrt{f(x)}-h\, \psi(x)\bigr)^2\,dx=o(h^2)$ and
$\int_{\RR}\bigl(\psi(x+h)-\psi(x)\bigr)^2\,dx=o(1)$, as well as that
$\bil\|\psi\bir\|>0$.

Note that under this assumptions, the model of i.i.d.\ observations with
density $f(x+\theta)$ is, in particular, LAN at $\theta=0$ with Fisher
information $I=4\,\bil\|\psi\bir\|^2=4\int_{\RR}\psi^2(x)\;dx$ $\bigl($see,
for example,~\cite[Chapter~2.1]{IKh}$\bigr)$ and so, using characteristic
functions, we have
$$
\lim_{n\to\infty}\Bigl(\eb
e^{it\ln\frac{f(\varepsilon+u/\sqrt{n})}{f(\varepsilon)}}\Bigr)^n
=e^{i\bigl(-\frac{Iu^2}2\bigr)t-\frac12 I u^2 t^2}
$$
and, more generally,
\begin{equation}
\label{assumption}
\lim_{\gamma\to0}\Bigl(\eb
e^{it\ln\frac{f(\varepsilon+\gamma)}{f(\varepsilon)}}\Bigr)^{1/\gamma^2}
=e^{i\bigl(-\frac I2\bigr)t-\frac12 I t^2}
\end{equation}
for all $t\in\RR$.

Note also, that only the convergence~\eqref{assumption} will be needed in our
considerations. So, one can rather assume it directly, or make any other
regularity assumptions sufficient for it as, for example, H\'ajek's conditions:
$f$ is differentiable and the Fisher information $I=\int_{\RR}
f^{-1}(x)\bigl(f'(x)\bigr)^2\,dx$ is finite and strictly positive $\bigl($see,
for example,~\cite[Chapter~2.2]{IKh}$\bigr)$.

Note finally, that in the Gaussian case the regularity assumptions clearly
hold and we have $I=1$.

\subsection*{The asymptotics \boldmath$\gamma\to 0$}

Let us consider the process $X_{\gamma,f}(y)=Z_{\gamma,f}(y/I\gamma^2)$,
$y\in\RR$, where $Z_{\gamma,f}$ is defined by~\eqref{proc1}. Note that
\begin{align*}
\frac{\int_{\RR}y\,X_{\gamma,f}(y)\;dy}{\int_{\RR}\,X_{\gamma,f}(y)\;dy}&=
I\gamma^2\zeta_{\gamma,f}\,,\\
\inf\Bigl\{z : X_{\gamma,f}(z)=\sup_{y\in\RR}X_{\gamma,f}(y)\Bigr\}&=I\gamma^2
\xi_{\gamma,f}^-\\
\intertext{and}
\sup\Bigl\{z : X_{\gamma,f}(z)=\sup_{y\in\RR}X_{\gamma,f}(y)\Bigr\}&=I\gamma^2
\xi_{\gamma,f}^+\,,
\end{align*}
where the random variables $\zeta_{\gamma,f}$ and $\xi_{\gamma,f}^{\pm}$ are
defined by~\eqref{vars1}. Remind also the process $Z_0$ on $\RR$ defined
by~\eqref{proc2} and the random variables $\zeta_0$ and $\xi_0$ defined
by~\eqref{vars2}. Recall finally the quantities
$B_{\gamma,f}=\eb\zeta_{\gamma,f}^2$,
$M_{\gamma,f}^\alpha=\eb(\xi_{\gamma,f}^\alpha)^2$,
$E_{\gamma,f}^\alpha=B_{\gamma,f}/M_{\gamma,f}^\alpha$, as well as
$B_0=\eb\zeta_0^2=16\,\zeta(3)$, $M_0=\eb\xi_0^2=26$ and
$E_0=B_0/M_0=8\,\zeta(3)/13$. Now we can state the main result of the present
paper.

\begin{theorem}
\label{T1}
The process $X_{\gamma,f}$ converges weakly in the space
$\mathcal{D}_0(-\infty ,+\infty)$ to the process $Z_0$ as $\gamma \to 0$. In
particular, the random variable $I\gamma^2\zeta_{\gamma,f}$ converges weakly
to the random variable $\zeta_0$ and, for any $\alpha\in[0,1]$, the random
variable $I\gamma^2\xi_{\gamma,f}^\alpha$ converges weakly to the random
variable $\xi_0$.  Moreover, for any $k>0$ we have
$$
I^k\gamma^{2k}\,\eb\zeta_{\gamma,f}^k \to \eb\zeta_0^k\quad\text{and}\quad
I^k\gamma^{2k}\,\eb(\xi_{\gamma,f}^\alpha)^k \to \eb\xi_0^k.
$$
In particular, $I^2\gamma^4 B_{\gamma,f} \to 16\,\zeta(3)$, $I^2\gamma^4
M_{\gamma,f}^\alpha \to 26$ and $E_{\gamma,f}^\alpha \to 8\,\zeta(3)/13$.
\end{theorem}

The results concerning the random variable $\zeta_{\gamma,f}$ are direct
consequence of~\cite[Theorem~1.10.2]{IKh} and the following three lemmas.

\begin{lemma}
\label{L1}
The finite-dimensional distributions of the process $X_{\gamma,f}$ converge to
those of $Z_0$ as $\gamma \to 0$.
\end{lemma}

\begin{lemma}
\label{L2}
For any $C>1/4$ we have
$$
\eb\left|X_{\gamma,f}^{1/2}(y_1)-X_{\gamma,f}^{1/2}(y_2)\right|^2\leq
C\abs{y_1-y_2}
$$
for all sufficiently small $\gamma$ and all $y_1,y_2\in\RR$.
\end{lemma}

\begin{lemma}
\label{L3}
For any $c\in\left]\,0\,{,}\;1/8\,\right[$ we have
$$
\eb X_{\gamma,f}^{1/2}(y)\leq\exp\bigl(-c\abs y\bigr)
$$
for all sufficiently small $\gamma$ and all $y\in\RR$.
\end{lemma}

Note that these lemmas are not sufficient to establish the weak convergence of
the process $X_{\gamma,f}$ in the space $\mathcal{D}_0(-\infty ,+\infty)$ and the
results concerning the random variable $\xi_{\gamma,f}^\alpha$. However, the
increments of the process $\ln X_{\gamma,f}$ being independent, the convergence of
its restrictions (and hence of those of $X_{\gamma,f}$) on finite intervals
$[A,B]\subset\RR$ $\bigl($that is, convergence in the Skorohod space
$\mathcal{D}[A,B]$ of functions on $[A,B]$ without discontinuities of the
second kind$\bigr)$ follows from~\cite[Theorem~6.5.5]{GS}, Lemma~\ref{L1} and
the following lemma.

\begin{lemma}
\label{L4}
For any $\delta>0$ we have
$$
\lim_{h\to 0}\ \lim_{\gamma\to 0}\ \sup_{\abs{y_1-y_2} < h}
\pb\Bigl\{\bigl|\ln X_{\gamma,f}(y_1)-\ln X_{\gamma,f}(y_2)\bigr|>\delta
\Bigr\}=0.
$$
\end{lemma}

Now, Theorem~\ref{T1} follows from the following estimate on the tails of the
process $X_{\gamma,f}$ by standard argument $\bigl($see, for
example,~\cite{IKh}$\bigr)$.

\begin{lemma}
\label{L5}
For any $b\in\left]\,0\,{,}\;1/12\,\right[$ we have
$$
\pb\biggl\{\sup_{\abs y>A} X_{\gamma,f}(y) > e^{-bA}\biggr\} \leq 4\,e^{-bA}
$$
for all sufficiently small $\gamma$ and all $A>0$.
\end{lemma}

The proofs of all these lemmas will be given in Section~\ref{PL}.

\subsection*{The asymptotics \boldmath$\gamma\to+\infty$}

Now let us discuss the second possible asymptotics $\gamma\to+\infty$. It can
be shown that in this case, the process $Z_{\gamma,f}$ converges weakly in the
space $\mathcal{D}_0(-\infty ,+\infty)$ to the process
$Z_\infty(x)=\1_{\{-\eta<x<\tau\}}$, $x\in\RR$, where $\eta$ and $\tau$ are
two independent exponential random variables with parameter $1$. So, the
random variables $\zeta_{\gamma,f}$, $\xi_{\gamma,f}^{-}$,
$\xi_{\gamma,f}^{+}$ and $\xi_{\gamma,f}^\alpha$ converge weakly to the random
variables
\begin{align*}
\zeta_\infty&=\frac{\int_{\RR}x\,Z_\infty(x)\;dx}{\int_{\RR}\,Z_\infty(x)\;dx}
=\frac{\tau-\eta}{2}\,,\\
\xi_\infty^-&=\inf\Bigl\{z : Z_\infty(z)=\sup_{x\in\RR}Z_\infty(x)\Bigr\}
=-\eta,\\
\xi_\infty^+&=\sup\Bigl\{z : Z_\infty(z)=\sup_{x\in\RR}Z_\infty(x)\Bigr\}
=\tau\\
\intertext{and}
\xi_\infty^\alpha&=\alpha\,\xi_\infty^-+(1-\alpha)\,\xi_\infty^+
=(1-\alpha)\,\tau -\alpha\,\eta,
\end{align*}
respectively. It can be equally shown that, moreover, for any $k>0$ we have
$$
\eb\zeta_{\gamma,f}^k\to\eb\zeta_\infty^k \quad\text{and}\quad
\eb(\xi_{\gamma,f}^\alpha)^k\to\eb(\xi_\infty^\alpha)^k.
$$
In particular, denoting $B_\infty=\eb\zeta_\infty^2$,
$M_\infty^\alpha=\eb(\xi_\infty^\alpha)^2$ and
$E_\infty^\alpha=B_\infty/M_\infty^\alpha$, we finally have
\begin{align}
B_{\gamma,f} &\to B_\infty=\eb\Bigl(\frac{\tau-\eta}{2}\Bigr)^2=\frac12\,,\notag\\
\label{b}
M_{\gamma,f}^\alpha &\to M_\infty^\alpha=\eb\bigl((1-\alpha)\,\tau
-\alpha\,\eta\bigr)^2=6\left(\alpha-\frac12\right)^2+\frac12\\
\intertext{and}
\label{c}
E_{\gamma,f}^\alpha &\to
E_\infty^\alpha=\frac{1}{12\left(\alpha-\frac12\right)^2+1}\,.
\end{align}

Let us note that these convergences are natural, since the process $Z_\infty$
can be considered as a particular case of the process $Z_{\gamma,f}$ with
$\gamma=+\infty$ under natural conventions $f(\varepsilon\pm\infty)=0$ and
$\ln0=-\infty$.

Note also, that $Z_\infty$ is the limiting likelihood ratio process in the
problem of estimating the parameter $\theta$ by i.i.d.\ uniform observations
on $[\theta,\theta+1]$. So, in this problem, the variables $\zeta_\infty$ and
$\xi_\infty^\alpha$ are the limiting distributions of the Bayesian estimators
and of the appropriately chosen maximum likelihood estimator, respectively,
while $B_\infty$ and $M_\infty^\alpha$ are the limiting mean squared errors of
these estimators and, the Bayesian estimators being asymptotically efficient,
$E_\infty^\alpha$ is the asymptotic relative efficiency of this maximum
likelihood estimator.

Finally observe, that the formulae~\eqref{b} and~\eqref{c} clearly imply that
in the latter problem (as well as in any problem having $Z_\infty$ as limiting
likelihood ratio) the best choice of the maximum likelihood estimator is
$\alpha=1/2$, and that the so chosen maximum likelihood estimator is
asymptotically efficient. This choice was also suggested for TAR model (which
has limiting likelihood ratio $Z_\gamma$) by Chan and Kutoyants
in~\cite{CK}. For large values of $\gamma$ this suggestion is confirmed by our
asymptotic results. However, we see that for small values of $\gamma$ the
choice of $\alpha$ will not be so important, since the limits in
Theorem~\ref{T1} do not depend on $\alpha$.

\subsection*{Numerical simulations}

Here we present some numerical simulations (in the Gaussian case) of the
quantities $B_\gamma$, $M_\gamma^\alpha$ and $E_\gamma^\alpha$ for
$\gamma\in\left]0,\infty\right[$. Besides giving approximate values of these
quantities, the simulation results illustrate both the asymptotics
$$
B_\gamma=\frac{B_0}{\gamma^4}+o(\gamma^{-4}),\quad
M_\gamma^\alpha=\frac{M_0}{\gamma^4}+o(\gamma^{-4})\quad\text{and}\quad
E_\gamma^\alpha\to E_0 \quad\text{as}\quad \gamma\to 0,
$$
with $B_0=16\,\zeta(3)\approx19.2329$, $M_0=26$ and
$E_0=8\,\zeta(3)/13\approx 0.7397$, and
$$
B_\gamma\to B_\infty,\quad M_\gamma^\alpha\to M_\infty^\alpha
\quad\text{and}\quad E_\gamma^\alpha\to E_\infty^\alpha \quad\text{as}\quad
\gamma\to\infty,
$$
with $B_\infty=0.5$, $M_\infty^\alpha=6\,(\alpha-0.5)^2+0.5$ and
$E_\infty^\alpha=1/\bigl(12\,(\alpha-0.5)^2+1\bigr)$.

First, we simulate the events $x_1^+,x_2^+,\ldots$ of the Poisson process
$\Pi_+$ and the events $x_1^-,x_2^-,\ldots$ of the Poisson process $\Pi_-$
$\bigl($both of intensity $1\bigr)$, as well as the partial sums
$S_1^+,S_2^+,\ldots$ of the i.i.d.\ $\cN(0,1)$ sequence
$\varepsilon_1^+,\varepsilon_2^+,\ldots$ and the partial sums
$S_1^-,S_2^-,\ldots$ of the i.i.d.\ $\cN(0,1)$ sequence
$\varepsilon_1^-,\varepsilon_2^-,\ldots.$ For convenience we also put
$x_0^+=x_0^-=S_0^+=S_0^-=0$.

Then we calculate
\begin{align*}
\zeta_\gamma&=\frac{\int_{\RR}x\,Z_\gamma(x)\;dx}
{\int_{\RR}\,Z_\gamma(x)\;dx}\\
&=\frac{\sum\limits_{i=0}^{\infty}\frac12\,e^{S_i^+}
\bigl({(x_{i+1}^+)}^2-{(x_i^+)}^2\bigr)-
\sum\limits_{i=0}^{\infty}\frac12\,e^{S_i^-}
\bigl({(x_{i+1}^-)}^2-{(x_i^-)}^2\bigr)}
{\sum\limits_{i=0}^{\infty}e^{S_i^+}(x_{i+1}^+-x_i^+)+
\sum\limits_{i=0}^{\infty}e^{S_i^-}(x_{i+1}^--x_i^-)}\;,\\
\xi_\gamma^-&=\inf\Bigl\{z :
Z_\gamma(z)=\sup_{x\in\RR}Z_\gamma(x)\Bigr\}=\begin{cases}
x_k^+, &\text{if } S_k^+ > S_\ell^-,\\
-x_{\ell+1}^-, &\text{otherwise},
\end{cases}\\
\xi_\gamma^+&=\sup\Bigl\{z :
Z_\gamma(z)=\sup_{x\in\RR}Z_\gamma(x)\Bigr\}=\begin{cases}
x_{k+1}^+, &\text{if } S_k^+ \geq S_\ell^-,\\
-x_\ell^-, &\text{otherwise},
\end{cases}\\
\intertext{and}
\xi_\gamma^\alpha&=\alpha\,\xi_\gamma^-+
(1-\alpha)\,\xi_\gamma^+,\\
\end{align*}
where
$$
k=\argmax_{i\geq 0}S_i^+\quad\text{and}\quad\ell=\argmax_{i\geq 0}S_i^-,
$$
and we use the values $1/2$, $1/4$ and $0$ for $\alpha$. Note that in this
Gaussian case (due to the symmetry of the process $Z_\gamma$) the random
variable $\xi_\gamma^{1-\alpha}$ has the same law as the variable
$-\xi_\gamma^\alpha$, that's why we use for $\alpha$ only values less or equal
than $1/2$.

Finally, repeating these simulations $10^7$ times (for each value of
$\gamma$), we approximate $B_\gamma=\eb\zeta_\gamma^2$ and
$M_\gamma^\alpha=\eb(\xi_\gamma^\alpha)^2$ by the empirical second moments,
and $E_\gamma^\alpha=B_\gamma/M_\gamma^\alpha$ by their ratio.

\begin{figure}[!h]
\centering\includegraphics*[width=0.6\textwidth]{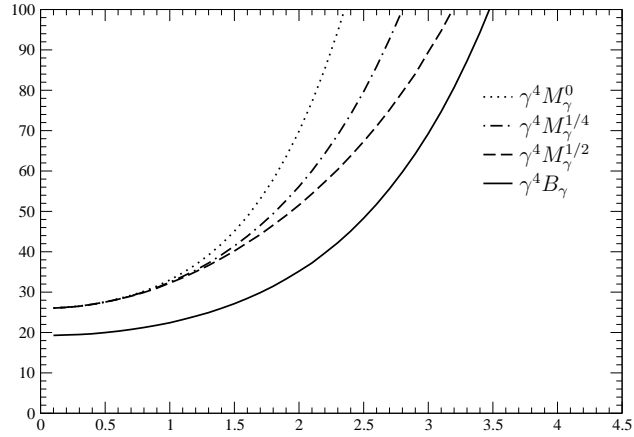}
\caption{$\gamma^4B_\gamma$ and $\gamma^4M_\gamma^\alpha$ ($\gamma\to 0$
  asymptotics)}
\label{fig1}
\end{figure}

\begin{figure}[!h]
\centering\includegraphics*[width=0.6\textwidth]{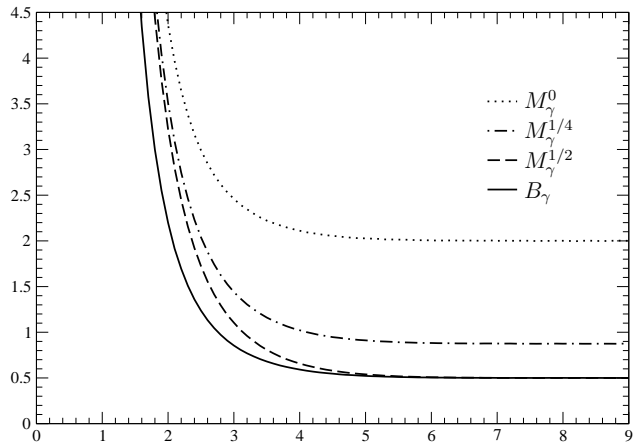}
\caption{$B_\gamma$ and $M_\gamma^\alpha$ ($\gamma\to\infty$ asymptotics)}
\label{fig2}
\end{figure}

\begin{figure}[!h]
\centering\includegraphics*[width=0.6\textwidth]{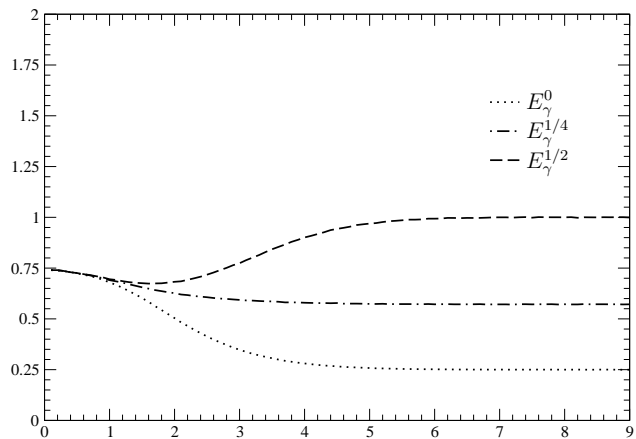}
\caption{$E_\gamma^\alpha$ (both asymptotics)}
\label{fig3}
\end{figure}

The results of the numerical simulations are presented in
Figures~\ref{fig1}--\ref{fig3}. The $\gamma\to 0$ asymptotics of the limiting
mean squared errors is illustrated in Figure~\ref{fig1}, where we rather
plotted the functions $\gamma^4 B_\gamma$ and $\gamma^4 M_\gamma^\alpha$,
making apparent the constants $B_0\approx19.2329$ and $M_0=26$. One can
observe here that the choice $\alpha=1/2$ is the best one, though its
advantage diminishes as~$\gamma$ approaches $0$ and seems negligible for
$\gamma<1$.

In Figure~\ref{fig2} we illustrate the $\gamma\to\infty$ asymptotics of the
limiting mean squared errors by plotting the functions $B_\gamma$ and
$M_\gamma^\alpha$ themselves. Here the advantage of the choice $\alpha=1/2$ is
obvious, and one can observe that for $\gamma>5$ this choice makes negligible
the loss of efficiency resulting from the use of the maximum likelihood
estimator instead of the asymptotically efficient Bayesian estimators.

Finally, in Figure~\ref{fig3} we illustrate the behavior both at $0$ and at
$\infty$ of the asymptotic relative efficiency of the maximum likelihood
estimators by plotting the functions $E_\gamma^\alpha$. All the observations
made above can be once more noticed in this figure. Note also that as $\gamma$
increases from $0$ to $\infty$, the asymptotic relative efficiency seems first
to decrease from $E_0\approx 0.7397$ for all the maximum likelihood
estimators, before increasing back to $E_\infty^\alpha$ for the maximum
likelihood estimators with $\alpha$ close to the optimal value $1/2$.

\section{Proofs of the lemmas}
\label{PL}

For the sake of clarity, for each lemma we will first give the proof in the
particular Gaussian case (in which it is more explicit) and then explain how
it can be extended to the general one.

\subsection*{Proof of Lemma~\ref{L1}}

Note that the restrictions of the process $\ln X_\gamma(y)=\ln
Z_\gamma(y/\gamma^2)$, $y\in\RR$, (as well as those of the process $\ln Z_0$)
on $\RR_+$ and on $\RR_-$ are mutually independent processes with stationary
and independent increments. So, to obtain the convergence of all the
finite-dimensional distributions, it is sufficient to show the convergence of
one-dimensional distributions only, that is, the weak convergence of $\ln
X_\gamma(y)$ to
$$
\ln Z_0(y)=W(y)-\frac{\abs{y}}{2}\sim\cN
\biggl(-\frac{\abs{y}}{2},\abs{y}\biggr)
$$
for all $y\in\RR$. Moreover, these processes being symmetric, it is sufficient
to consider $y\in\RR_+$ only.

The characteristic function $\varphi_\gamma(t)$ of $\ln X_\gamma(y)$ is
$$
\begin{aligned}
\varphi_\gamma(t)&=\eb\,e^{it\ln X_\gamma(y)}=\eb\,e^{-it\gamma
\sum_{k=1}^{\Pi_+\bigl(y/\gamma^2\bigr)}\varepsilon_k^+
-it\frac{\gamma^2}{2}\Pi_+(y/\gamma^2)}\\
&=\eb\,\eb\,\Bigl(e^{-it\gamma\sum_{k=1}^{\Pi_+\bigl(y/\gamma^2\bigr)}
\varepsilon_k^+-it\frac{\gamma^2}{2}\Pi_+(y/\gamma^2)}\Bigm|\sF_{\Pi_+}\Bigr)\\
&=\eb\biggl(e^{-it\frac{\gamma^2}{2}\Pi_+(y/\gamma^2)}
\prod_{k=1}^{\Pi_+(y/\gamma^2)}\eb\,e^{-it\gamma\varepsilon_k^+}\biggr)\\
&=\eb\,e^{-it\frac{\gamma^2}{2}\Pi_+(y/\gamma^2)-\frac{t^2\gamma^2}{2}
\Pi_+(y/\gamma^2)}=\eb\,e^{-\frac{\gamma^2}{2}(it+t^2)\Pi_+(y/\gamma^2)}
\end{aligned}
$$
where we have denoted $\sF_{\Pi_+}$ the $\sigma$-algebra related to the
Poisson process $\Pi_+$, used the independence of $\varepsilon_k^+$ and
$\Pi_+$ and recalled that $\eb\,e^{it\varepsilon}=e^{-t^2/2}$.

Then, noting that $\Pi_+(y/\gamma^2)$ is a Poisson random variable of
parameter $y/\gamma^2$ with moment generating function
$\eb\,e^{t\Pi_+(y/\gamma^2)}=\exp\bigl(\frac{y}{\gamma^2}(e^t-1)\bigr)$, we
get
$$
\begin{aligned}
\ln\varphi_\gamma(t)&=\frac{y}{\gamma^2}
\Bigl(e^{-\frac{\gamma^2}{2}(it+t^2)}-1\Bigr)=\frac{y}{\gamma^2}
\Bigl(-\frac{\gamma^2}{2}(it+t^2)+o(\gamma^2)\Bigr)\\
&=-\frac{y}{2}(it+t^2)+o(1)\to -\frac{y}{2}(it+t^2)=\ln\eb\,e^{it\ln Z_0(y)}
\end{aligned}
$$
as $\gamma\to0$ and so, in the Gaussian case Lemma~\ref{L1} is proved.

\bigskip

In the general case, proceeding similarly we get
$$
\begin{aligned}
\varphi_\gamma(t)&=\eb\,e^{it\ln X_{\gamma,f}(y)}=
\eb\,e^{it\sum_{k=1}^{\Pi_+\bigl(y/I\gamma^2\bigr)}\ln
\frac{f(\varepsilon_k^++\gamma)}{f(\varepsilon_k^+)}}\\
&=\eb\biggl(\Bigl(\eb\,e^{it\ln\frac{f(\varepsilon+\gamma)}
{f(\varepsilon)}}\Bigr)^{\Pi_+(y/I\gamma^2)}\biggr)\to
e^{i\bigl(-\frac y2\bigr)t-\frac12 y t^2}=\eb\,e^{it\ln Z_0(y)}\\
\end{aligned}
$$
by dominated convergence theorem, since
$$
\Bigl(\eb\,e^{it\ln\frac{f(\varepsilon+\gamma)}
{f(\varepsilon)}}\Bigr)^{1/\gamma^2}\to
e^{i\bigl(-\frac I2\bigr)t-\frac12 I t^2}
$$
by~\eqref{assumption}, and $\gamma^2\,\Pi_+(y/I\gamma^2)$ converges clearly to
$y/I$ in $L^2$ (and hence in probability).

\subsection*{Proof of Lemma~\ref{L3}}

Now we turn to the proof of Lemma~\ref{L3} (we will prove Lemma~\ref{L2} just
after). For $y>0$ we have
$$
\begin{aligned}
\eb X_\gamma^{1/2}(y)&=\eb\,\eb\,\Bigl(e^{-\frac\gamma2
\sum_{k=1}^{\Pi_+\bigl(y/\gamma^2\bigr)}\varepsilon_k^+
-\frac{\gamma^2}{4}\Pi_+(y/\gamma^2)}\Bigm|\sF_{\Pi_+}\Bigr)\\
&=\eb\,e^{-\frac{\gamma^2}{4}\Pi_+(y/\gamma^2)+\frac{\gamma^2}{8}
\Pi_+(y/\gamma^2)}=\eb\,e^{-\frac{\gamma^2}{8}\Pi_+(y/\gamma^2)}\\
&=\exp\biggl(\frac{y}{\gamma^2}\Bigl(e^{-\frac{\gamma^2}8}-1\Bigr)\biggr).
\end{aligned}
$$

The process $X_\gamma$ being symmetric, we have
\begin{equation}
\label{Eroot}
\eb X_\gamma^{1/2}(y)=\exp\biggl(\frac{\abs{y}}{\gamma^2}
\Bigl(e^{-\frac{\gamma^2}8}-1\Bigr)\biggr)
\end{equation}
for all $y\in\RR$ and, since
$$
\frac{1}{\gamma^2}\Bigl(e^{-\frac{\gamma^2}8}-1\Bigr)=
\frac{1}{\gamma^2}\Bigl(-\frac{\gamma^2}8+o(\gamma^2)\Bigr)\to-\frac18
$$
as $\gamma\to0$, for any $c\in\left]\,0\,{,}\;1/8\,\right[$ we have $\eb
X_\gamma^{1/2}(y)\leq\exp\bigl(-c\abs{y}\bigr)$ for all sufficiently small
$\gamma$ and all $y\in\RR$. So, in the Gaussian case Lemma~\ref{L3} is proved.

\bigskip

In the general case, equality~\eqref{Eroot} becomes $\eb
X_{\gamma,f}^{1/2}(y)=\exp\bigl(\abs y(I_\gamma-1)/I\gamma^2\bigr)$ with
$$
I_\gamma=\eb\sqrt{\frac{f(\varepsilon+\gamma)}{f(\varepsilon)}}
\leq\sqrt{\eb\,\frac{f(\varepsilon+\gamma)}{f(\varepsilon)}}=1.
$$
Recall the convergence~\eqref{assumption} of characteristic functions and note
that $I_\gamma^{1/\gamma^2}$ are the corresponding moment generating functions
at point $1/2$. The convergence of these moment generating functions (at any
point smaller than~$1$) follows from the fact that for all $\gamma$ they are
equal $1$ at point $1$ (which provides uniform integrability). Thus we have
$I_\gamma^{1/\gamma^2}\to e^{-I/8}$, which implies $(\ln
I_\gamma)/\gamma^2\to-I/8$, and so $(I_\gamma-1)/I\gamma^2\to-1/8$.

\subsection*{Proof of Lemma~\ref{L2}}

First we consider the case $y_1,y_2\in\RR_+$ (say $y_1\geq
y_2$). Using~\eqref{Eroot} and taking into account the stationarity and the
independence of the increments of the process $\ln X_\gamma$ on $\RR_+$, we
can write
$$
\begin{aligned}
\eb\left|X_\gamma^{1/2}(y_1)-X_\gamma^{1/2}(y_2)\right|^2&=\eb X_\gamma(y_1)+\eb
X_\gamma(y_2)-2\, \eb X_\gamma^{1/2}(y_1)X_\gamma^{1/2}(y_2)\\
&=2-2\,\eb X_\gamma(y_2)\,\eb\,\frac{X_\gamma^{1/2}(y_1)}{X_\gamma^{1/2}(y_2)}\\
&=2-2\,\eb X_\gamma^{1/2}\bigl(\abs{y_1-y_2}\bigr)\\
&=2-2\exp\biggl(\frac{\abs{y_1-y_2}}{\gamma^2}
\Bigl(e^{-\frac{\gamma^2}8}-1\Bigr)\biggr)\\
&\leq-2\,\frac{\abs{y_1-y_2}}{\gamma^2}\Bigl(e^{-\frac{\gamma^2}8}-1\Bigr)\leq
\frac14\abs{y_1-y_2}.
\end{aligned}
$$

The process $X_\gamma$ being symmetric, we have the same result for the case
$y_1,y_2\in\RR_-$.

Finally, if $y_1y_2\leq 0$ (say $y_2\leq 0\leq y_1$), we have
$$
\begin{aligned}
\eb\left|X_\gamma^{1/2}(y_1)-X_\gamma^{1/2}(y_2)\right|^2&=2-2\,\eb
X_\gamma^{1/2}(y_1)\,\eb X_\gamma^{1/2}(y_2)\\
&=2-2\exp\biggl(\frac{\abs{y_1}}{\gamma^2}\Bigl(e^{-\frac{\gamma^2}8}-1\Bigr)
+\frac{\abs{y_2}}{\gamma^2}\Bigl(e^{-\frac{\gamma^2}8}-1\Bigr)\biggr)\\
&=2-2\exp\biggl(\frac{\abs{y_1-y_2}}{\gamma^2}
\Bigl(e^{-\frac{\gamma^2}8}-1\Bigr)\biggr)\\
&\leq\frac14\abs{y_1-y_2},
\end{aligned}
$$
and so, in the Gaussian case we obtain even more than the assertion of
Lemma~\ref{L2}.

\bigskip

In the general case, proceeding similarly we get
$$
\eb\left|X_{\gamma,f}^{1/2}(y_1)-X_{\gamma,f}^{1/2}(y_2)\right|^2\leq
-2\,\frac{\abs{y_1-y_2}}{I\gamma^2}(I_\gamma-1)
$$
and, since $-2(I_\gamma-1)/I\gamma^2\to1/4$, the proof is concluded.

\subsection*{Proof of Lemma~\ref{L4}}

First let $y_1,y_2\in\RR_+$ (say $y_1\geq y_2$) such that
$\Delta=\abs{y_1-y_2}<h$. Then, noting that conditionally to $\sF_{\Pi_+}$ the
random variable
$$
\ln X_\gamma(\Delta)=-\gamma
\sum_{k=1}^{\Pi_+(\Delta/\gamma^2)}\varepsilon_k^+
-\frac{\gamma^2}{2}\Pi_+(\Delta/\gamma^2)
$$
is Gaussian with mean
$-\frac{\gamma^2}{2}\Pi_+(\Delta/\gamma^2)$ and variance
$\gamma^2\Pi_+(\Delta/\gamma^2)$, we get
$$
\begin{aligned}
\pb\Bigl\{\bigl|\ln X_\gamma(y_1)-\ln X_\gamma(y_2)\bigr|>\delta\Bigr\}
&\leq\frac1{\delta^2}\,\eb\bigl|\ln X_\gamma(y_1)-\ln X_\gamma(y_2)\bigr|^2\\
&=\frac1{\delta^2}\,\eb\bigl|\ln X_\gamma(\Delta)\bigr|^2\\
&=\frac1{\delta^2}\,\eb\,\eb\,\Bigl(\bigl(\ln X_\gamma(\Delta)\bigr)^2
\Bigm|\sF_{\Pi_+}\Bigr)\\
&=\frac1{\delta^2}\,\eb\left(\gamma^2\Pi_+(\Delta/\gamma^2)
+\frac{\gamma^4}{4}\bigl(\Pi_+(\Delta/\gamma^2)\bigr)^2\right)\\
&=\frac1{\delta^2}\,\Biggl(\Delta+\frac{\gamma^4}{4}\biggl(
\frac{\Delta}{\gamma^2}+\frac{\Delta^2}{\gamma^4}\biggr)\Biggr)\\
&=\frac1{\delta^2}\,\bigl((1+\gamma^2/4)\Delta+\Delta^2/4\bigr)\\
&<\frac1{\delta^2}\,\bigl(\beta(\gamma)\,h+h^2/4\bigr)
\end{aligned}
$$
where $\beta(\gamma)=1+\gamma^2/4\to 1$ as $\gamma\to 0$. So, we have
$$
\begin{aligned}
\lim_{\gamma\to 0}\ \sup_{\abs{y_1-y_2} < h}\pb\Bigl\{\bigl|\ln X_\gamma(y_1)-\ln
X_\gamma(y_2)\bigr|>\delta\Bigr\} &\leq\lim_{\gamma\to 0}
\frac1{\delta^2}\,\bigl(\beta(\gamma)\,h+h^2/4\bigr)\\
&=\frac1{\delta^2}\left(h+\frac{h^2}4\right),
\end{aligned}
$$
and hence
$$
\lim_{h\to 0}\ \lim_{\gamma\to 0}\ \sup_{\abs{y_1-y_2} < h}\pb\Bigl\{\bigl|\ln
X_\gamma(y_1)-\ln X_\gamma(y_2)\bigr|>\delta\Bigr\}=0,
$$
where the supremum is taken only over $y_1,y_2\in\RR_+$.

The process $X_\gamma$ being symmetric, we have the same conclusion with the
supremum taken over $y_1,y_2\in\RR_-$.

Finally, for $y_1y_2\leq 0$ (say $y_2\leq 0\leq y_1$) such that
$\abs{y_1-y_2}<h$, using the elementary inequality $(a-b)^2\leq 2(a^2+b^2)$ we
get
$$
\begin{aligned}
\pb\Bigl\{\bigl|\ln X_\gamma(y_1)-\ln X_\gamma(y_2)\bigr|>\delta\Bigr\}
&\leq\frac1{\delta^2}\,\eb\bigl|\ln X_\gamma(y_1)-\ln X_\gamma(y_2)\bigr|^2\\
&\leq\frac2{\delta^2}\biggl(\eb\bigl|\ln X_\gamma(y_1)\bigr|^2+
\eb\Bigl|\ln X_\gamma\bigl(\abs{y_2}\bigr)\Bigr|^2\biggr)\\
&=\frac2{\delta^2}\bigl(\beta(\gamma)y_1+y_1^2/4
+\beta(\gamma)\abs{y_2}+\abs{y_2}^2/4\bigr)\\
&<\frac2{\delta^2}\Bigl(\beta(\gamma)h+h^2/4\Bigr),
\end{aligned}
$$
which again yields the desired conclusion. So, in the Gaussian case
Lemma~\ref{L4} is proved.

\bigskip

Another way to prove this lemma, is to notice first that the weak convergence
of $\ln X_\gamma(y)$ to $\ln Z_0(y)$ (established in Lemma~\ref{L1}) is
uniform with respect to $y\in K$ for any compact $K\subset\RR$. Indeed, the
uniformity of the convergence of the characteristic functions in the proof of
Lemma~\ref{L1} is obvious, and so one can apply, for example, Theorem~7 from
Appendix~I of~\cite{IKh}, whose remaining conditions are easily checked.

Second, using this uniformity we obtain
$$
\begin{aligned}
\lim_{\gamma\to 0}\ \sup_{\abs{y_1-y_2} < h}
\pb\Bigl\{\bigl|\ln X_\gamma(y_1)-\ln X_\gamma(y_2)\bigr|>\delta\Bigr\}
&=\lim_{\gamma\to 0}\ \sup_{\abs{y} < h}
\pb\Bigl\{\bigl|\ln X_\gamma(y)\bigr|>\delta\Bigr\}\\
&=\sup_{\abs{y} < h}\pb\Bigl\{\bigl|\ln Z_0(y)\bigr|>\delta\Bigr\}
\end{aligned}
$$
where the supremum is taken over $y_1,y_2\in\RR$ such that $y_1y_2\geq 0$, and
$$
\lim_{\gamma\to 0}\ \sup_{\abs{y_1-y_2} < h}
\pb\Bigl\{\bigl|\ln X_\gamma(y_1)-\ln X_\gamma(y_2)\bigr|>\delta\Bigr\}
\leq2\sup_{\abs{y} < h}\pb\Bigl\{\bigl|\ln Z_0(y)\bigr|>\frac\delta2\Bigr\}
$$
where the supremum is taken over $y_1,y_2\in\RR$ such that $y_1y_2\leq 0$.

Finally, reminding that $\ln Z_0(y)\sim\cN\bigl(-\abs{y}/2\,,\,\abs{y}\bigr)$
and denoting $\Phi$ the distribution function of the standard Gaussian law, we
get
$$
\begin{aligned}
\pb\Bigl\{\bigl|\ln Z_0(y)\bigr|>\delta\Bigr\}
&=\Phi\biggl(-\frac\delta{\sqrt{\abs y}}+\frac{\sqrt{\abs y}}2\,\biggr)+1-
\Phi\biggl(\frac\delta{\sqrt{\abs y}}+\frac{\sqrt{\abs y}}2\,\biggr)\\
&\leq\Phi\biggl(-\frac\delta{\sqrt{h}}+\frac{\sqrt{h}}2\,\biggr)+1-
\Phi\biggl(\frac\delta{\sqrt{h}}\biggr)
\end{aligned}
$$
for $\abs{y}<h$. The last expression does not depend on $y$ and clearly
converges to $0$ as $h\to 0$, so the assertion of the lemma follows.

It remains to observe that this second proof does not use any particularity of
the process $X_\gamma$ and, hence, is trivially extendable to the general
case.

\subsection*{Proof of Lemma~\ref{L5}}

Taking into account the symmetry of the process $\ln X_\gamma$, as well as the
stationarity and the independence of its increments on $\RR_+$, we obtain
\begin{equation}
\label{Psup}
\begin{aligned}
\pb\biggl\{\sup_{\abs y>A} X_\gamma(y) > e^{-bA}\biggr\}&\leq
2\,\pb\biggl\{\sup_{y>A} X_\gamma(y) > e^{-bA}\biggr\}\\
&\leq 2\,e^{\,bA/2}\;
\eb\sup_{y>A} X_\gamma^{1/2}(y)\\
&=2\,e^{\,bA/2}\;\eb X_\gamma^{1/2}(A)\;
\eb\sup_{y>A}\frac{X_\gamma^{1/2}(y)}{X_\gamma^{1/2}(A)}\\
&=2\,e^{\,bA/2}\;\eb X_\gamma^{1/2}(A)\;\eb\sup_{z>0} X_\gamma^{1/2}(z).
\end{aligned}
\end{equation}

In order to estimate the last factor we write
$$
\begin{aligned}
\eb\sup_{z>0} X_\gamma^{1/2}(z)&=\eb\exp\left(\frac12\,
\sup_{z>0}\,\Biggl(-\gamma\sum_{k=1}^{\Pi_+(z/\gamma^2)}\varepsilon_k^+
-\frac{\gamma^2}{2}\Pi_+(z/\gamma^2)\Biggr)\right)\\
&=\eb\exp\left(\frac12\,
\sup_{n\in\NN}\,\biggl(-\gamma\sum_{k=1}^{n}\varepsilon_k^+
-\frac{n\gamma^2}{2}\biggr)\right).
\end{aligned}
$$
Now, let us observe that the random walk $S_n=-\sum_{k=1}^{n}\varepsilon_k^+$,
$n\in\NN$, has the same law as the restriction on $\NN$ of a standard Brownian
motion $W$. So,
$$
\begin{aligned}
\eb\sup_{z>0} X_\gamma^{1/2}(z)&=\eb\exp\left(\frac12\,
\sup_{n\in\NN}\,\bigl(\gamma W(n)-n\gamma^2/2\bigr)\right)\\
&=\eb\exp\left(\frac12\,\sup_{n\in\NN}\,\bigl(W(n\gamma^2)
-n\gamma^2/2\bigr)\right)\\
&\leq\eb\exp\left(\frac12\,\sup_{t>0}\,\bigl(W(t)-t/2\bigr)\right)=
\eb\exp\left(\frac12\,M\right)
\end{aligned}
$$
with an evident notation. It is known that the random variable $M$ is
exponential of parameter $1$ $\bigl($see, for example,~\cite{BS}$\bigr)$ and
hence, using its moment generating function $\eb\,e^{tM}=(1-t)^{-1}$, we get
\begin{equation}
\label{Esup}
\eb\sup_{z>0} X_\gamma^{1/2}(z) \leq 2.
\end{equation}

Finally, taking $b\in\left]\,0\,{,}\;1/12\,\right[$ we have
$3b/2\in\left]\,0\,{,}\;1/8\,\right[$ and, combining~\eqref{Psup},
\eqref{Esup} and using Lemma~\ref{L3}, we finally obtain
$$
\begin{aligned}
\pb\biggl\{\sup_{\abs y>A} X_\gamma(y) > e^{-bA}\biggr\}&\leq
4\,e^{\,bA/2}\,\exp\Bigl(-\frac{3b}{2}A\Bigr)=4\,e^{-bA}
\end{aligned}
$$
for all sufficiently small $\gamma$ and all $A>0$, which concludes the proof
in the Gaussian case.

\bigskip

In the general case the proof is almost the same. Note that we have no longer
the symmetry of the process $X_{\gamma,f}$, so we need to consider the cases
$y>A$ and $y<-A$ separately. Besides that, the only difference is in the
derivation of the bound~\eqref{Esup}. Here we get
$$
\eb\sup_{z>0} X_{\gamma,f}^{1/2}(z)=\eb\exp\left(\frac12\,M\right),
$$
where $M$ is the supremum of the random walk $S_n=\sum_{k=1}^{n}X_k$,
$n\in\NN$, with $X_k=\ln\frac{f(\varepsilon_k^++\gamma)}
{f(\varepsilon_k^+)}\,$. Note that
$$
\eb\,e^{X_1}=\eb\,\frac{f(\varepsilon+\gamma)} {f(\varepsilon)}=1,
$$
and so, the cummulant generating function $k(t)=\ln(\eb\,e^{tX_1})$ of $X_1$
admits a strictly positive zero $t_0=1$. Hence, by the well-known
Cram\'er-Lundberg bound on the tail probabilities of $M$ $\bigl($see, for
example, Theorem~5.1 from Chapter~XIII of~\cite{Asm}$\bigr)$, we have
$$
\pb(M>x)\leq e^{-t_0\,x}=e^{-x}
$$
for all $x>0$. Finally, denoting $F$ the distribution function of $M$ and
using this bound we obtain
$$
\begin{aligned}
\eb\exp\left(\frac12\,M\right)&=\int_{\RR}e^{\,x/2}\;dF(x)\\
&=\Bigl[e^{\,x/2}\bigl(F(x)-1\bigr)\Bigr]_{-\infty}^{+\infty}-\;
\frac12\int_{\RR}e^{\,x/2}\bigl(F(x)-1\bigr)\;dx\\
&=\frac12\int_{\RR_-}e^{\,x/2}\;dx+
\frac12\int_{\RR_+}e^{\,x/2}\bigl(1-F(x)\bigr)\;dx\\
&\leq 1+\frac12\int_{\RR_+}e^{-x/2}\;dx=2,
\end{aligned}
$$
which concludes the proof.


\end{document}